\documentclass[12pt]{article}

\usepackage{amsmath}
\usepackage{amssymb}
\usepackage[dvips]{graphicx}

\begin{document}
\title{{\bf 
{\large Amphicheiral links with special properties, I}}
\footnotetext[0]{%
2010 {\it Mathematics Subject Classification}: 
57M25, 57M27.\\
{\it Keywords}:
algebraically split link; amphicheiral link; invertible link; Alexander polynomial; 
Reidemeister torsion.
}}
\author{{\footnotesize 
Teruhisa KADOKAMI}}
\date{{\footnotesize July 10, 2011}\bigskip\\
}
\maketitle

\newcommand{\circlenum}[1]{{\ooalign{%
\hfill$\scriptstyle#1$\hfill\crcr$\bigcirc$}}}

\newcommand{\svline}[1]{\multicolumn{1}{|c}{#1}}
\newfont{\bg}{cmr10 scaled\magstep4}
\newcommand{\bigzerol}{\smash{\hbox{\bg 0}}}
\newcommand{\bigzerou}{\smash{\lower1.7ex\hbox{\bg 0}}}

\newcommand{\bsquare}{\hbox{\rule{6pt}{6pt}}}
\newcommand{\qed}{\hbox{\rule[-2pt]{3pt}{6pt}}}
\newcommand{\Int}{\mathrm{Int}\ \! }
\newcommand{\Ker}{\mathrm{Ker}\ \! }
\newcommand{\Ig}{\mathrm{Im}\ \! }
\newcommand{\aug}{\mathrm{aug}\ \!}
\newcommand{\pj}{\mathrm{pr}\ \!}
\newcommand{\Tor}{\mathrm{Tor}\ \!}
\newcommand{\Spin}{\mathrm{Spin}\ \!}
\newcommand{\Eul}{\mathrm{Eul}\ \!}
\newcommand{\Vect}{\mathrm{Vect}\ \!}
\newcommand{\HULL}{\mathrm{HULL}\ \!}
\newcommand{\real}{\mathrm{real}\ \!}
\newcommand{\rank}{\mathrm{rank}\ \!}
\newcommand{\ord}{\mathrm{ord}\ \!}
\newcommand{\Sign}{\mathrm{Sign}\ \!}
\newcommand{\Hom}{\mathrm{Hom}\ \!}
\newcommand{\ad}{\mathrm{ad}\ \!}
\newcommand{\Det}{\mathrm{Det}\ \!}
\newcommand{\lk}{\mathrm{lk}\ \!}
\newcommand{\pt}{\mathrm{pt}}
\newcommand{\al}{$\alpha$}
\newcommand{\dis}{\displaystyle}

\newtheorem{df}{Definition}[section]
\newtheorem{lm}[df]{Lemma}
\newtheorem{theo}[df]{Theorem}
\newtheorem{re}[df]{Remark}
\newtheorem{pr}[df]{Proposition}
\newtheorem{ex}[df]{Example}
\newtheorem{co}[df]{Corollary}
\newtheorem{cl}[df]{Claim}
\newtheorem{qu}[df]{Question}
\newtheorem{pb}[df]{Problem}
\newtheorem{cj}[df]{Conjecture}

\makeatletter
\renewcommand{\theequation}{%
\thesection.\arabic{equation}}
\@addtoreset{equation}{section}
\makeatother

\begin{abstract}
{\footnotesize 
\setlength{\baselineskip}{10pt}
\setlength{\oddsidemargin}{0.25in}
\setlength{\evensidemargin}{0.25in}
We provide necessary conditions for the Alexander polynomials
of algebraically split component-preservingly amphicheiral links.
We raise a conjecture that the Alexander polynomial of an algebraically split 
component-preservingly amphicheiral link with even components is zero.
Our necessary conditions and some examples support the conjecture.}
\end{abstract}

\section{Introduction}\label{sec:intro}
Let $L=K_1\cup \cdots \cup K_r$ be an oriented 
$r$-component link in $S^3$ with $r \ge 1$.
For an oriented knot $K$, we denote the orientation-reversed knot by $-K$.
If $\varphi$ is an orientation-reversing (orientation-preserving, respectively)
homeomorphism of $S^3$ so that
$\varphi(K_i)=\varepsilon_{\sigma(i)} K_{\sigma(i)}$ for all $i=1, \ldots, r$
where $\varepsilon_i=+$ or $-$,
and $\sigma$ is a permutation of $\{1, 2, \ldots, r\}$, 
then $L$ is said 
an {\it $(\varepsilon_1, \ldots, \varepsilon_r; \sigma)$-amphicheiral link}
(an {\it $(\varepsilon_1, \ldots, \varepsilon_r; \sigma)$-invertible link},
respectively).
A term ``amphicheiral link"  is used as a general term for
an $(\varepsilon_1, \ldots, \varepsilon_r; \sigma)$-amphicheiral link.
A link is said an {\it interchangeable link} if it is
an $(\varepsilon_1, \ldots, \varepsilon_r; \sigma)$-invertible link
such that $\sigma$ is not the identity.
An $(\varepsilon_1, \ldots, \varepsilon_r; \sigma)$-invertible link
is said an {\it invertible link} simply if there exists $1\le i\le r$
such that $\varepsilon_i=-$.
If $\sigma$ is the identity, then an amphicheiral link is said
a {\it component-preservingly amphicheiral link}, and
$\sigma$ may be omitted from the notation.
We mainly deal with component-preservingly amphicheiral link
in the present paper.
If every $\varepsilon_i=\varepsilon$ is identical for all $i=1, \ldots, r$
(including the case that $\sigma$ is not the identity), then 
an $(\varepsilon_1, \ldots, \varepsilon_r; \sigma)$-amphicheiral link
(an $(\varepsilon_1, \ldots, \varepsilon_r; \sigma)$-invertible link,
respectively)
is said an $(\varepsilon)$-amphicheiral link
(an $(\varepsilon)$-invertible link, respectively).
We use the notations $+=+1=1$ and $-=-1$.
A link $L$ with at least $2$-component is said 
an {\it algebraically split link} if the linking number of
every $2$-component sublink of $L$ is zero.
We note that 
a component-preservingly $(\varepsilon)$-amphicheiral link
is an algebraically split link.

\medskip

Necessary conditions for the Alexander polynomials of amphicheiral knots 
are studied by R.~Hartley \cite{Ha}, R.~Hartley and A.~Kawauchi \cite{HK},
and A.~Kawauchi \cite{Kw2}.
In \cite{Kw2}, non-invertibility of $8_{17}$ is firstly proved by the conditions.
On the other hand, T.~Sakai \cite{Sa} proved that
any one-variable Laurent polynomial $f(t)$ over $\mathbb{Z}$
such that $f(t)=f(t^{-1})$ and $f(1)=1$ is realized by
the Alexander polynomial of a strongly invertible knot in $S^3$.
B.~Jiang, X.~Lin, Shicheng Wang and Y.~Wu \cite{JLWW}
showed that 
(1) a twisted Whitehead doubled knot is amphicheiral
if and only if it is the unknot or the figure eight knot, and
(2) a prime link with at least $2$ components and up to $9$ crossings is 
component-preservingly $(+)$-amphicheiral if and only if it is the Borromean rings.
They used S.~Kojima and M.~Yamasaki's $\eta$-function \cite{KY}.
Shida Wang \cite{Wa} determined 
prime component-preservingly $(+)$-amphicheiral links 
with at least $2$ components and up to $11$ crossings by the same method 
as \cite{JLWW}. 
There are four such links.
For geometric studies of symmetries of arborescent knots, see
F.~Bonahon and L.~C.~Siebenmann \cite{BS}.
In \cite{Kd1}, we determined symmetries such as invertibility,
amphicheirality and interchangeability of $2$-bridge links
by the parameters such as Schubert's normal form,
Conway's normal form and Conway's normal form
whose entries are even integers.

\medskip

In the present paper, we study necessary conditions for 
the Alexander polynomials of algebraically split 
component-preservingly amphicheiral links 
by computing the Reidemeister torsions of surgered manifolds along the link.
The results and the techniques of the present paper
have been already applied to some directions.
The author and A.~Kawauchi \cite{KK} obtained a necessary condition
by invariants deduced from the quadratic form of a link \cite{Kw1, Kw5},
showed a partial affirmative answer for the conjecture stated below,
and determined prime amphicheiral links with up to $9$ crossings
by combining with the present results.
The author \cite{Kd2} developed methods to detect component-preservingly
amphicheiral links by the results in \cite{JLWW} and \cite{Ha, HK, Kw2},
and determined prime amphicheiral links with up to $11$ crossings
by techniques including the present results.
There are 27 prime amphicheiral links with up to $11$ crossings.
The techniques of the present paper are based on V.~G.~Turaev \cite{Tr}.
By the same techniques, the following two results on Dehn surgeries
are shown:
In \cite{KMS}, the author, N.~Maruyama and M.~Shimozawa determined
all Dehn surgeries yieding lens spaces (i.e.\ lens surgeries).
In \cite{Kd3},
we will show that the $\lambda$-component Milnor link 
with $\lambda \ge 4$ does not have a lens surgery
by the Reidemeister torsion and some geometric techniques.

\medskip

Let ${\mit \Delta}_L={\mit \Delta}_L(t_1, \ldots, t_r)$ 
be the Alexander polynomial of $L$ which is an element of 
an $r$-variable Laurent polynomial ring 
${\mit \Lambda}_r:=
\mathbb{Z}[t_1^{\pm 1}, \ldots, t_r^{\pm 1}]$
over $\mathbb{Z}$ where $t_i$\ ($i=1, 2, \ldots, r$)
is a variable corresponding to a meridian of $K_i$.

\medskip

We raise a conjecture:

\begin{cj}\label{cj:cj1}
For an even-component algebraically split 
component-preservingly amphicheiral link $L$,
we have ${\mit \Delta}_L=0$.
\end{cj}

Our results stated below support the conjecture.
In Section \ref{sec:rem}, we explain about other supporting results 
in \cite{Kd2, KK}.
We remark that the similar statement for the odd-component case does not hold.
For example, a link which is connected sums of 
copies of the Borromean rings is a component-preservingly 
amphicheiral link with odd components,
and has the non-zero Alexander polynomial.
For any odd number, there are such examples.
For an algebraically split component-preservingly amphicheiral link $L$,
if $L'$ is obtained from $L$ by taking untwisted parallels of some components,
and the number of components of $L'$ is strictly greater than that of $L$,
then $L'$ is also an algebraically split component-preservingly amphicheiral link, and
we have ${\mit \Delta}_{L'}\equiv 0$.
Therefore we cannot find a counterexample for the conjecture
by the construction.
We raise supporting examples in Section \ref{sec:rem}.

\medskip

If $L$ is an $r$-component algebraically split link with $r \ge 2$, 
then the Alexander polynomial of $L$ is of the form:
\begin{equation*}\label{eq:L}
{\mit \Delta}_L={\mit \Delta}_L(t_1, \ldots, t_r)\doteq
(t_1-1)\cdots (t_r-1)f(t_1, \ldots, t_r)
\end{equation*}
where we can take
$f=f(t_1, \ldots, t_r)\in {\mit \Lambda}_r$
satisfying 
$f(t_1, \ldots, t_r)=f(t_1^{-1}, \ldots, t_r^{-1})$.
Note that $f(t_1, \ldots, t_r)$ is uniquely determined up to
multiplication of $\pm 1$.
We set $I_r=\{1, 2, \ldots, r\}$.
If we take $I=\{i_1, i_2, \ldots, i_s\}\subset I_r$, then
we denote $L_I=K_{i_1}\cup \cdots \cup K_{i_s}$, $|I|=s$,
and ${\mit \Lambda}_I
=\mathbb{Z}[t_{i_1}^{\pm 1}, \ldots, t_{i_s}^{\pm 1}]
\cong {\mit \Lambda}_s$.
If $s \ge 2$, then 
\begin{equation*}\label{eq:LI}
{\mit \Delta}_{L_I}
={\mit \Delta}_{L_I}(t_{i_1}, \ldots, t_{i_s})\doteq
\prod_{i\in I}(t_i-1) f_I(t_{i_1}, \ldots, t_{i_s})
\end{equation*}
where we can take
$f_I=f_I(t_{i_1}, \ldots, t_{i_s})\in {\mit \Lambda}_I$
satisfying 
$f_I(t_{i_1}, \ldots, t_{i_s})=f_I(t_{i_1}^{-1}, \ldots, t_{i_s}^{-1})$,
and the sign of $f_I$ is uniquely determined by 
(\ref{eq:Torres1}) and Lemma \ref{lm:Torres2}.
In particular, $f=f_{I_r}$.
We set $u(I)=(u_i)_{i\in I_r\setminus I}$ where $u_i\in \{1, -1\}$.
For $J=\{j_1, \ldots, j_k\} \supset I$, 
we use the following notations:
$$\begin{array}{cl}
k_J(u(I)) & : \mbox{the number of $1$ in $u_i\ (i\in J\setminus I)$},
\medskip\\
\eta_J(u(I)) & =(-1)^{k_J(u(I))},
\medskip\\
F_J(I) & : \mbox{the polynomial obtained by substituting
$t_i=1$\ $(i\in J\setminus I)$ to $f_J$},
\medskip\\
F(I) & =F_{I_r}(I).
\end{array}$$
We set
$$S_{u(I)}^{\mathrm{even}}
=\sum_{
{\scriptstyle J\supset I, |J\setminus I|\mathrm{: even}}
\atop
{\scriptstyle 2\le |J| \le r}}
\eta_J(u(I))F_J(I),$$
and 
$$S_{u(I)}^{\mathrm{odd}}
=\sum_{
{\scriptstyle J\supset I, |J\setminus I|\mathrm{: odd}}
\atop
{\scriptstyle 2\le |J| \le r}}
\eta_J(u(I))F_J(I).$$
The following is our first main theorem:

\begin{theo}\label{th:MT1}
Let $L=K_1\cup \cdots \cup K_r$ be
an $r$-component algebraically split component-preservingly amphicheiral link
where $r \ge 2$, and
$I\subset I_r$.
Then for any $u(I)$, we have the following:
\begin{enumerate}
\item[(1)]
If $|I|=1$, then $S_{u(I)}^{\mathrm{odd}}=0$.

\item[(2)]
If $2\le |I|\le r-1$, then
$S_{u(I)}^{\mathrm{even}}=0$ or
$S_{u(I)}^{\mathrm{odd}}=0$.

\end{enumerate}
\end{theo}

To prove the theorem,
we compute the Reidemeister torsions of the manifolds 
surgered along $L$ with associated coefficients to $u(I)$.
We can deduce some corollaries.

\begin{co}\label{co:F=0}
Under the same assumption as in Theorem \ref{th:MT1},
we have the following:
\begin{enumerate}
\item[(1)]
If $r$ is even and $|I|=1$, then $F(I)=0$.

\item[(2)]
If $I=I_r\setminus \{i\}$ (i.e.\ $|I|=r-1$)
and ${\mit \Delta}_{L_I}\ne 0$, 
then $f$ is divisible by $t_i-1$.

\item[(3)]
If $I=I_r\setminus \{i\}$ (i.e.\ $|I|=r-1$)
and $F(I)\ne 0$, 
then ${\mit \Delta}_{L_I}=0$.

\end{enumerate}
\end{co}

In particular, if $r=2$, then we have the following:

\begin{co}\label{co:2-comp1}
If $L=K_1\cup K_2$ is
an algebraically split component-preservingly amphicheiral link, 
then ${\mit \Delta}_L$ is divisible by
$(t_1-1)^2(t_2-1)^2$.
\end{co}

The following is our second main theorem:

\begin{theo}\label{th:MT2}
If $L=K_1\cup \cdots \cup K_r$ is
an $r$-component algebraically split 
component-preservingly
$(\varepsilon)$-amphicheiral link
with $r$ even, and $\varepsilon=+$ or $-$, 
then the Alexander polynomial of $L$ satisfies
${\mit \Delta}_L(t^{\eta_1}, \ldots, t^{\eta_r})=0$
where $\eta_i\in \{1, -1\}$\ $(i=1, \ldots, r)$.
\end{theo}

To prove the theorem, we span a Seifert surface on $L$.
This method is a slightly extended argument
in \cite[Theorem 2.1]{Ha}.
In particular, if $r=2$, then we have the following:

\begin{co}\label{co:2-comp2}
If $L=K_1\cup K_2$ is
an algebraically split $(\varepsilon, \varepsilon)$-amphicheiral link
where $\varepsilon=+$ or $-$,
then ${\mit \Delta}_L$ is divisible by
$(t_1-1)^2(t_2-1)^2(t_1t_2-1)(t_1-t_2)$.
\end{co}

We remark that after proving the results above,
the author and A.~Kawauchi \cite{KK} showed
that Conjecture \ref{cj:cj1} is afffirmative
for even-component algebraically split
component-preservingly $(\varepsilon)$-amphicheiral links
by invariants deduced from the quadratic form of a link \cite{Kw1, Kw5}.

\medskip

In Section \ref{sec:R-tor}, we prepare facts on the Alexander polynomials.
In Section \ref{sec:ampinv}, we discuss about basic properties
on amphicheiral links and invertible links.
We give an almost purely algebraic proof for a lemma
due to Hartley \cite{Ha}.
In Section \ref{sec:MT1}, we prove Theorem \ref{th:MT1}
and its corollaries.
In Section \ref{sec:MT2}, we prove Theorem \ref{th:MT2}.
In Section \ref{sec:rem}, 
we raise some examples which support Conjecture \ref{cj:cj1}.

\section{Alexander polynomials as Reidemeister torsions}\label{sec:R-tor}
Let $X$ be a finite CW complex, and
$\psi : \mathbb{Z}[H_1(X)]\to R$ is a ring homomorphism
where $R$ is an integral domain.
Then we denote the Reidemeister torsion of $X$
related to $\psi$ by $\tau^{\psi}(X)\in Q(R)$ 
where $Q(R)$ is the quotient field of $R$ (see \cite{Tr}).
The value $\tau^{\psi}(X)$ is determined up to multiplication
of $\pm \psi(h)\ (h\in H_1(X))$.
An equation between two values $A$ and $B\in Q(R)$ 
is denoted by $A\doteq B$
if $A=\pm \psi(h)B$ for some $h\in H_1(X)$.
When $\psi$ is the identity, we denote $\tau^{\psi}(X)$ by $\tau(X)$.

\medskip

The Alexander polynomial is a kind of the Reidemeister torsion.

\begin{lm}\label{lm:Alexander}{\rm (\cite{Mi3, Tr})}
Let $L=K_1\cup \cdots \cup K_r$ be an $r$-component link,
and $E_L$ the complement of $L$.
Then we have
$$\tau(E_L)\doteq \left\{
\begin{array}{cl}
{\displaystyle \frac{{\mit \Delta}_L(t_1)}{t_1-1}} & (r =1),\medskip\\
{\mit \Delta}_L(t_1, \ldots, t_r) & (r \ge 2).
\end{array}
\right.$$
\end{lm}

We will use the surgery formula to show Theorem 1.1.

\begin{lm}\label{lm:surgery}{\rm (surgery formula)}
Let $M_0$ be a compact $3$-manifold whose boundary consists
of tori, $V$ a solid torus whose core is $l'$, and
$M=M_0\cup_f V$ is the result of Dehn filling
where $f : \partial V\to \partial M_0$ is an attaching map.
Let $\psi : \mathbb{Z}[H_1(M)]\to R$ be a ring homomorphism
where $R$ is an integral domain,
and $\psi_0 : \mathbb{Z}[H_1(M_0)]\to R$ the induced map from $\psi$.
If $\psi([l'])\ne 1$, then we have
$$\tau^{\psi}(M)\doteq
\tau^{\psi_0}(M_0)(\psi([l'])-1)^{-1}.$$
\end{lm}

We raise some properties on the Alexander polynomials.

\begin{lm}\label{lm:dual}{\rm (duality \cite{Kw3, Mi3, To, Tr})}
Let $L=K_1\cup \cdots \cup K_r$ be an $r$-component link.
Then we have
$${\mit \Delta}_L(t_1)=t_1^a{\mit \Delta}_L(t_1^{-1})\quad
(r=1)$$
where $a$ is even, and
$${\mit \Delta}_L(t_1, \ldots, t_r)=
(-1)^rt_1^{a_1}\cdots t_r^{a_r}
{\mit \Delta}_L(t_1^{-1}, \ldots, t_r^{-1})
\quad (r \ge 2)$$
where $a_i\equiv 1+\sum_{j\ne i}\mathrm{lk}\ \! (K_i, K_j)\ 
(\mathrm{mod}\ \! 2)$\ $(i=1, \ldots, r)$.
\end{lm}

The Torres condition is a special case of the surgery formula.

\begin{lm}\label{lm:Torres}{\rm (Torres condition \cite{Kw3, To, Tr})}
Let $L=K_1\cup \cdots \cup K_r\cup K_{r+1}$ 
be an oriented $(r+1)$-compnent link, and
$L'=K_1\cup \cdots \cup K_r$ an $r$-component sublink.
Then we have
$${\mit \Delta}_L(t_1, 1)\doteq
\frac{t_1^{\ell}-1}{t_1-1}
{\mit \Delta}_{L'}(t_1)\quad (r=1)$$
where $\ell$ is the linking number of $K_1$ and $K_2$, and
$${\mit \Delta}_L(t_1, \ldots, t_r, 1)\doteq
(t_1^{\ell_1}\cdots t_r^{\ell_r}-1)
{\mit \Delta}_{L'}(t_1, \ldots, t_r)\quad
(r \ge 2)$$
where $\ell_i$ is the linking number of $K_i$ and $K_{r+1}$\ 
$(i=1, \ldots, r)$.
\end{lm}

One necessary condition for the Alexander polynomial of
an amphicheiral or invertible link is the following:

\begin{lm}\label{lm:amp}
Let $L=K_1\cup \cdots \cup K_r$ be
an $r$-component 
$(\varepsilon_1, \ldots, \varepsilon_r)$-amphicheiral or invertible link
where $\varepsilon_i=+$ or $-$\ $(i=1, \ldots, r)$.
Then we have
$${\mit \Delta}_L(t_1, \ldots, t_r)\doteq 
{\mit \Delta}_L(t_1^{\varepsilon_1}, \ldots, 
t_r^{\varepsilon_r}).$$
\end{lm}

\section{Amphicheiral link and invertible link}\label{sec:ampinv}
We raise basic properties of amphicheiral links and invertible links.

\begin{lm}\label{lm:sub}
Let $L=K_1\cup \cdots \cup K_r$ be an $r$-component link.
\begin{enumerate}
\item[(1)]
If $L$ is an $(\varepsilon_1, \ldots, \varepsilon_r)$-amphicheiral link, then
a sublink $L'=K_{i_1}\cup \cdots \cup K_{i_s}$\ 
$(1\le i_1<\cdots<i_s\le r)$ is an
$(\varepsilon_{i_1}, \ldots, \varepsilon_{i_s})$-amphicheiral link.

\item[(2)]
If $L'=K_{i_1}\cup \cdots \cup K_{i_s}$\ 
$(1\le i_1<\cdots <i_s\le r)$
is an $s$-component sublink of $L$ such that $s \ge 3$ is odd, and
$\ell_{1,2}\cdot \ell_{2,3}\cdots \ell_{s-1,s}\cdot \ell_{s,1}\ne 0$
where $\ell_{p,q}$ is the linking number of $K_{i_p}$ and $K_{i_q}$,
then $L$ is not component-preservingly amphicheiral.

\item[(3)]
If $L$ is an $(\varepsilon_1, \ldots, \varepsilon_r)$-invertible link, then
a sublink $L'=K_{i_1}\cup \cdots \cup K_{i_s}$\ 
$(1\le i_1<\cdots<i_s\le r)$ is an
$(\varepsilon_{i_1}, \ldots, \varepsilon_{i_s})$-invertible link.

\item[(4)]
Let $L=K_1\cup K_2$ be a $2$-component link with non-zero linking number.
If $L$ is an invertible link, then $L$ is a $(-, -)$-invertible link.

\end{enumerate}
\end{lm}

The linking numbers are the first information to detect 
both amphicheirality and invertibility as in Lemma \ref{lm:sub} (2) and (4).
By Lemma \ref{lm:sub} (1) and (3), to study sublinks
is also important for the problems.
Hartley \cite{Ha} showed the following by 
the JSJ (Jaco-Shalen-Johanson) decomposition.
The result is also about relation between amphicheirality and the linking number.
We reprove it by another way.

\begin{lm}\label{lm:even}{\rm (Hartley \cite{Ha})}
Let $L=K_1\cup K_2$ be a $2$-component link
with non-zero even linking number.
Then $L$ is not component-preservingly amphicheiral.
\end{lm}

\noindent
{\bf Proof}\ Suppose that $L$ is component-preservingly amphicheiral, and
the linking number of $K_1$ and $K_2$ is non-zero and even.
By Lemma \ref{lm:dual}, we may assume
\begin{equation}\label{eq:dual}
{\mit \Delta}_L(t_1, t_2)=
t_1t_2{\mit \Delta}_L(t_1^{-1}, t_2^{-1})
\end{equation}
By Lemma \ref{lm:amp} and Lemma \ref{lm:sub} (4), we may assume
$${\mit \Delta}_L(t_1, t_2)=
\eta t_1^{b_1}t_2^{b_2}{\mit \Delta}_L(t_1^{-1}, t_2)$$
where $\eta=+$ or $-$, and $b_1, b_2\in \mathbb{Z}$.
By substituting $t_2=1$ to (\ref{eq:dual}), 
we have $\eta=+$ and $b_1=1$.
By substituting $t_1=1$ to (\ref{eq:dual}), 
we have $b_2=0$, and hence we have
\begin{equation*}
{\mit \Delta}_L(t_1, t_2)=
t_1{\mit \Delta}_L(t_1^{-1}, t_2)
\end{equation*}
By substituting $t_1=-1$ to the equation above, 
we have ${\mit \Delta}_L(-1, t_2)=0$.
In the similar way, we have ${\mit \Delta}_L(t_1, -1)=0$, and hence
${\mit \Delta}_L(t_1, t_2)$ is divisible by $(t_1+1)(t_2+1)$.
We set ${\mit \Delta}_L(t_1, t_2)=(t_1+1)(t_2+1)g(t_1, t_2)$
where $g(t_1, t_2)\in {\mit \Lambda}_2$.

\medskip

By substituting $t_2=1$ to ${\mit \Delta}_L(t_1, t_2)$,
and Lemma \ref{lm:Torres}, we have
$${\mit \Delta}_L(t_1, 1)=2(t_1-1)g(t_1, 1)
\doteq \frac{t_1^{\ell}-1}{t_1-1}{\mit \Delta}_{K_1}(t_1).$$
Since the righthand side is not divisible by $2$,
we have a contradiction.
\qed

\bigskip

Both Lemma \ref{lm:sub} and Lemma \ref{lm:even}
motivate us to study amphicheirality of algebraically split links.
We remark that our proof of Lemma \ref{lm:even} works for the case
that $L$ is in an integral homology sphere 
with an orientation-reversing autohomeomorphism.

\section{Proof of Theorem \ref{th:MT1},
Corollary \ref{co:F=0} and Corollary \ref{co:2-comp1}}\label{sec:MT1}
To prove Theorem \ref{th:MT1}, we study the form of the Alexander polynomial
of an algebraically split link, and compute the Reidemeister torsions
of surgered manifolds along the link.

\medskip

Let $L=K_1\cup \cdots \cup K_r$ be an oriented
$r$-component algebraically split link where $r \ge 2$.
We add one component $K_{r+1}$ to $L$
such that $L_i=K_i\cup K_{r+1}$\ $(i=1, \ldots, r)$
is the connected sum of $K_i$ and the Hopf link,
where the linking number of $K_i$ and $K_{r+1}$ is $1$.
Then we have
\begin{equation}\label{eq:Alex1}
{\mit \Delta}_{L_i}(t_i, t_{r+1})
\doteq {\mit \Delta}_{K_i}(t_i)
\end{equation}
We set $\overline{L}=L\cup K_{r+1}$.
By Lemma \ref{lm:Torres}, we have
\begin{equation}\label{eq:Torres1}
\begin{matrix}
{\mit \Delta}_{\overline{L}}(t_1, \ldots, t_{r+1})
& = & (t_1\cdots t_r-1)(t_1-1)\cdots (t_r-1)
f(t_1, \ldots, t_r)\medskip\\
& & +(t_{r+1}-1)g(t_1, \ldots, t_{r+1})\hfill
\end{matrix}
\end{equation}
where $f(t_1, \ldots, t_r)\in {\mit \Lambda}_r$ and
$g(t_1, \ldots, t_{r+1})\in {\mit \Lambda}_{r+1}$.
By Lemma \ref{lm:dual}, we may assume that
\begin{equation}\label{eq:dual1}
f(t_1, \ldots, t_r)=f(t_1^{-1}, \ldots, t_r^{-1})
\end{equation}
and
\begin{equation}\label{eq:dual2}
{\mit \Delta}_{\overline{L}}(t_1, \ldots, t_{r+1})
=(-1)^{r+1}t_1^2\cdots t_r^2t_{r+1}^a
{\mit \Delta}_{\overline{L}}(t_1^{-1}, \ldots, t_{r+1}^{-1})
\end{equation}
where $a\equiv r+1\ (\mathrm{mod}\ \! 2)$.
For $I=\{i_1, i_2, \ldots, i_{\mu}\}\subset I_r$,
we set $L_I=K_{i_1}\cup \cdots \cup K_{i_s}$,
$\overline{L}_I=L_I\cup K_{r+1}$, and
$g_I\in {\mit \Lambda}_{\overline{I}}$ is obtained by substituting 
$t_j=1$ for all $j\in I_r\setminus I$
to $g(t_1, \ldots, t_{r+1})$.

\medskip

By (\ref{eq:Alex1}) and (\ref{eq:dual2}),
if $I=\{i\}$\ $(s=1)$, then we may take
\begin{equation}\label{eq:Torres2}
{\mit \Delta}_{\overline{L}_I}(t_i, t_{r+1})
={\mit \Delta}_{K_i}(t_i)
\end{equation}
where ${\mit \Delta}_{K_i}(t_i)=t_i^2{\mit \Delta}_{K_i}(t_i^{-1})$.
If $2\le s \le r$\ $(r \ge 2)$, then we may take
\begin{equation}\label{eq:Torres3}
\begin{matrix}
{\mit \Delta}_{\overline{L}_I}(t_{i_1}, \ldots, t_{i_s}, t_{r+1})
& = & {\displaystyle \left( \prod_{i\in I}t_i-1\right)\prod_{i\in I}(t_i-1)
f_I(t_{i_1}, \ldots, t_{i_s})}\hfill \medskip\\
&  & {\displaystyle +(t_{r+1}-1)
g_I'(t_{i_1}, \ldots, t_{i_s}, t_{r+1})}\hfill
\end{matrix}
\end{equation}
where
\begin{equation*}
{\mit \Delta}_{L_I}(t_{i_1}, \ldots, t_{i_s})
=\prod_{i\in I}(t_i-1)f_I(t_{i_1}, \ldots, t_{i_s}),
\end{equation*}
\begin{equation*}
f_I(t_{i_1}, \ldots, t_{i_s})=f_I(t_{i_1}^{-1}, \ldots, t_{i_s}^{-1})
\in {\mit \Lambda}_I
\end{equation*}
and
\begin{equation*}
g_I'(t_{i_1}, \ldots, t_{i_s}, t_{r+1})
\in {\mit \Lambda}_{\overline{I}}.
\end{equation*}
We set $f_I=f_I(t_{i_1}, \ldots, t_{i_s})$ and
$g_I'=g_I'(t_{i_1}, \ldots, t_{i_s}, t_{r+1})$.
We remark that 
$f_{I_r}=f(t_1, \ldots, t_r)$ and
$g_{I_r}'=g(t_1, \ldots, t_r)$.

\begin{lm}\label{lm:Torres2}
Under the situation above, for $1\le s \le r-1$, we have
\begin{equation*}
g_I=(t_{r+1}-1)^{r-s-1}
{\mit \Delta}_{\overline{L}_I}(t_{i_1}, \ldots, t_{i_s}, t_{r+1}).
\end{equation*}
\end{lm}

\noindent
{\bf Proof}\ 
By applying Lemma \ref{lm:Torres} repeatedly, we have the result.
\qed

\bigskip

We exapand $g(t_1, \ldots, t_{r+1})$-part
in (\ref{eq:Torres1}) as follows:

\begin{lm}\label{lm:form}
If $r \ge 2$, then we have
\begin{equation*}
\begin{matrix}
g(t_1, \ldots, t_{r+1})
& = & {\displaystyle (t_{r+1}-1)^{r-2}
\prod_{i=1}^r{\mit \Delta}_{K_i}(t_i)}\hfill \medskip\\
& & {\displaystyle +\sum_{
{\scriptstyle I\subset I_r}
\atop
{\scriptstyle 2\le s=|I| \le r-1}}
\prod_{i\in I}(t_i-1)
(t_{r+1}-1)^{r-s-1}}\hfill\medskip\\
& & {\displaystyle 
\cdot \left\{\left( \prod_{i\in I}t_i-1\right)f_I
+(t_{r+1}-1)h_I\right\}}\hfill\medskip\\
& & {\displaystyle +\prod_{i=1}^r(t_i-1)h}\hfill
\end{matrix}
\end{equation*}
where $f_I\in {\mit \Lambda}_I$,
$h_I\in {\mit \Lambda}_{\overline{I}}$,
and $h\in {\mit \Lambda}_{r+1}$.

\end{lm}

\noindent
{\bf Proof}\ 
We show by induction on $r$.

\medskip

\noindent
(i)\ The case $r=2$.

\medskip

By Lemma \ref{lm:Torres2} and (\ref{eq:Torres2}),
$g(t_1, t_2, t_3)-{\mit \Delta}_{K_1}(t_1){\mit \Delta}_{K_2}(t_2)$
is divisible by $(t_1-1)(t_2-1)$.
Hence we have the result.

\medskip

\noindent
(ii)\ The case $r \ge 3$.

\medskip

Suppose the case $r-1$.
By (\ref{eq:Torres3}), Lemma \ref{lm:Torres2} and the assumption, we have
\begin{equation*}
\begin{matrix}
g(t_1, \ldots, t_{r+1})
& = & {\displaystyle 
{\mit \Delta}_{\overline{L}_I}(t_1, \ldots, t_{r-1}, t_{r+1})
+(t_r-1)H_I}\hfill \medskip\\
& = & {\displaystyle 
\left(\prod_{i\in I}t_i-1\right) \prod_{i\in I}(t_i-1)f_I
+(t_{r+1}-1)g_I'
+(t_r-1)H_I}\hfill\medskip\\
& = & {\displaystyle \left(\prod_{i\in I}t_i-1\right) \prod_{i\in I}(t_i-1)f_I
+(t_{r+1}-1)^{r-2}\prod_{i\in I}{\mit \Delta}_{K_i}(t_i)
}\hfill\medskip\\
& & {\displaystyle +\sum_{
{\scriptstyle J\subset I}
\atop
{\scriptstyle 2\le k=|J| \le r-2}}
\prod_{i\in J}(t_i-1)
(t_{r+1}-1)^{r-k-1}}\hfill\medskip\\
& & {\displaystyle \cdot \left\{
\left( \prod_{i\in J}t_i-1\right)f_J
+(t_{r+1}-1)h_J\right\}}\hfill
\medskip\\
&  & {\displaystyle +\prod_{i\in I}(t_i-1)h_I
+(t_r-1)H_I}\hfill
\end{matrix}
\end{equation*}
where $I=I_r\setminus \{r \}$, 
$f_I\in {\mit \Lambda}_I$,
$g_I'\in {\mit \Lambda}_{\overline{I}}$,
$h_J\in {\mit \Lambda}_{\overline{J}}$,
and $H_I\in {\mit \Lambda}_{r+1}$.

\medskip

We set
\begin{equation*}
\begin{matrix}
G(t_1, \ldots, t_{r+1})
& = & {\displaystyle 
g(t_1, \ldots, t_{r+1})-(t_{r+1}-1)^{r-2}
\prod_{i=1}^r{\mit \Delta}_{K_i}(t_i)}\hfill \medskip\\
& & {\displaystyle -\sum_{
{\scriptstyle I\subset I_r}
\atop
{\scriptstyle 2\le s=|I| \le r-1}}
\prod_{i\in I}(t_i-1)
(t_{r+1}-1)^{r-s-1}}\hfill\medskip\\
& & {\displaystyle 
\cdot \left\{\left( \prod_{i\in I}t_i-1\right)f_I
+(t_{r+1}-1)h_I\right\}.}\hfill
\end{matrix}
\end{equation*}
Then we have
$G(t_1, \ldots, t_{r-1}, 1, t_{r+1})=0.$
In the similar way, if we substitute $t_i=1$ for any $1\le i\le r$
to $G(t_1, \ldots, t_{r+1})$, then we have $0$,
and hence $G(t_1, \ldots, t_{r+1})$ is divisible by 
$\prod_{i=1}^r(t_i-1)$.
\qed

\bigskip

Let $E$ be the complement of $\overline{L}$, and
$M=(\overline{L}; u_1, \ldots, u_{r+1})$
the result of $(u_1, \ldots, u_{r+1})$-surgery along $\overline{L}$
where $u_i\in \mathbb{Q}\cup \{\infty, \emptyset\}$.
For $I=\{i_1, i_2, \ldots, i_s\}\subset I_r$,
we suppose that $u_i=\emptyset$ if $i\in I$, 
$u_i\in \mathbb{Q}\cup \{\infty\}$ if $i\in I_r\setminus I$,
and $u_{r+1}\in \{\infty, \emptyset\}$.
We set $u(I)=(u_i)_{i\in I_r\setminus I}$ and
$u(\overline{I})=u(I)\cup \{u_{r+1}\}$.
If $u_{r+1}=\emptyset$, then we set $M=M_{u(\overline{I})}$, and
if $u_{r+1}=\infty$, then we set $M=M_{u(I)}$.
We set the natural inclusion 
$\iota : M_{u(\overline{I})}\hookrightarrow M_{u(I)}$.
From now on,
we consider only the case $u_i\in \{1, -1\}$\ $(i\in I_r\setminus I)$.

\medskip

Let $m_i$ and $l_i$\ $(i=1, \ldots, r+1)$ be 
a meridian and a longitude of $K_i$.
We denote the homology class of a loop $\gamma$ by $[\gamma]$.
In $H_1(E)$, we have
\begin{equation}\label{eq:relation1}
\begin{matrix}
[l_i]=[m_{r+1}]\hfill & (i=1, \ldots, r)\hfill\medskip\\
[l_{r+1}]=[m_1]\cdots [m_r]\hfill & 
\end{matrix}
\end{equation}
In $H_1(M_{u(\overline{I})})$, we have
\begin{equation}\label{eq:relation2}
[m_i']=[m_i]^{r_i}[l_i]=[m_i]^{r_i}[m_{r+1}]=1,\quad
[l_i']=[m_i]\qquad (i\in I_r\setminus I)
\end{equation}
where $m_i'$ and $l_i'$ are a meridian and a longitude
of the attaching solid torus for $K_i$\ $(i\in I_r\setminus I)$.
In $H_1(M_{u(I)})$, we have
\begin{equation}\label{eq:relation3}
[m_{r+1}']=[m_{r+1}]=1,\quad
[l_{r+1}']=[l_{r+1}]
\end{equation}
We note that
$$H_1(E)=\langle t_1, \ldots, t_r, t_{r+1}\rangle
\cong \mathbb{Z}^{r+1},$$
$$H_1(M_{u(\overline{I})})
=\langle t_{i_1}, \ldots, t_{i_s}, t_{r+1}\rangle
\cong \mathbb{Z}^{s+1},$$
and
$$H_1(M_{u(I)})=\langle t_{i_1}, \ldots, t_{i_s}\rangle
\cong \mathbb{Z}^s,$$
where we set $t_i=[m_i]$\ $(i=1, \ldots, r+1)$.

\medskip

For $J=\{j_1, \ldots, j_k\}\supset I$
and $p\in {\mit \Lambda}_{\overline{J}}$,
we use the following notations:
$$\begin{array}{cl}
k_J(u(I)) & : \mbox{the number of $1$ in $u_i\ (i\in J\setminus I)$},
\medskip\\
\rho_J(u(I)) & =(-1)^{k_J(u(I))}t_{r+1}^{-k_J(u(I))}, 
\medskip\\
\eta_J(u(I)) & =(-1)^{k_J(u(I))},
\medskip\\
\sigma_J(u(I)) & : \mbox{the sum of $u_i$\ $(i\in J\setminus I)$},
\medskip\\
p(u(I)) & \mbox{the polynomial obtained by substituting
$t_i=t_{r+1}^{-u_i}$\ $(i\in J\setminus I)$ to $p$}, 
\medskip\\
F_J(I) & : \mbox{the polynomial obtained by substituting
$t_i=1$\ $(i\in J\setminus I)$ to $f_J$},
\medskip\\
F(I) & =F_{I_r}(I).
\end{array}$$

The Reidemeister torsions of 
$M_{u(\overline{I})}$ is the following:

\begin{lm}\label{lm:R-tor1}
Suppose that $r \ge 2$ and $1\le |I|=s \le r-1$.
\begin{enumerate}
\item[(1)]
If $I=\{x\}$\ $(s=1)$, then we have:
\begin{equation*}
\begin{matrix}
\tau(M_{u(\overline{I})}) & \doteq &
{\displaystyle 
{\mit \Delta}_{K_x}(t_x)
\prod_{i\in I_r\setminus \{x\}}
{\mit \Delta}_{K_i}(t_{r+1}^{-u_i})}\hfill\medskip\\
& & {\displaystyle +(t_x-1)\sum_{
{\scriptstyle J\ni x}
\atop
{\scriptstyle 2\le |J| \le r}}
\rho_J(u(I))\left( t_xt_{r+1}^{-\sigma_J(u(I))}-1\right)
f_J(u(I))}\hfill\medskip\\
& & {\displaystyle 
+(t_{r+1}-1)Q_I}\hfill
\end{matrix}
\end{equation*}
where $Q_I\in {\mit \Lambda}_{\overline{I}}$.

\item[(2)]
If $2\le |I| \le r-1$, then we have:
\begin{equation*}
\begin{matrix}
\tau(M_{u(\overline{I})}) & \doteq &
{\displaystyle \prod_{i\in I}(t_i-1)
\sum_{
{\scriptstyle J\supset I}
\atop
{\scriptstyle 2\le |J| \le r}}
\rho_J(u(I))\left( \prod_{i\in I}t_it_{r+1}^{-\sigma_J(u(I))}-1\right)
f_J(u(I))}\hfill\medskip\\
& & {\displaystyle 
+(t_{r+1}-1)Q_I}\hfill
\end{matrix}
\end{equation*}
where $Q_I\in {\mit \Lambda}_{\overline{I}}$.

\end{enumerate}
\end{lm}

\noindent
{\bf Proof}\ 
By Lemma \ref{lm:Alexander}, Lemma \ref{lm:surgery}
and (\ref{eq:relation2}), we have
\begin{equation}\label{eq:tor1}
\begin{matrix}
\tau(M_{u(\overline{I})}) & \doteq &
{\displaystyle
{\mit \Delta}_{\overline{L}}(u(I))
\prod_{i\in I_r\setminus I}(t_{r+1}^{-u_i}-1)^{-1}}\hfill
\medskip\\
& \doteq & 
{\displaystyle
{\mit \Delta}_{\overline{L}}(u(I))
(t_{r+1}-1)^{-(r-s)}}\hfill
\end{matrix}
\end{equation}
By (\ref{eq:Torres1}), we have
\begin{equation}\label{eq:tor2}
\begin{matrix}
{\mit \Delta}_{\overline{L}}(u(I)) & \doteq &
{\displaystyle
\left(
\prod_{i\in I}t_it_{r+1}^{-\sigma_J(u(I))}-1
\right)
\prod_{i\in I}(t_i-1)
\prod_{i\in I_r\setminus I}
(t_{r+1}^{-u_i}-1)}f(u(I))\hfill \medskip\\
&  & 
{\displaystyle
+(t_{r+1}-1)g(u(I))}
\hfill \medskip\\
& = &
{\displaystyle
\rho_J(u(I))
\left(
\prod_{i\in I}t_it_{r+1}^{-\sigma_J(u(I))}-1
\right)
\prod_{i\in I}(t_i-1)
(t_{r+1}-1)^{r-s}}f(u(I))\hfill \medskip\\
&  & 
{\displaystyle
+(t_{r+1}-1)g(u(I))}
\hfill\\
\end{matrix}
\end{equation}
By Lemma \ref{lm:form}, we have
\begin{equation}\label{eq:tor3}
\begin{matrix}
g(u(I))
& = & {\displaystyle (t_{r+1}-1)^{r-2}
\prod_{i\in I}{\mit \Delta}_{K_i}(t_i)
\prod_{i\in I_r \setminus I}
{\mit \Delta}_{K_i}(t_{r+1}^{-u_i})}
\hfill \medskip\\
& & {\displaystyle +\sum_{
{\scriptstyle J\subset I_r}
\atop
{\scriptstyle 2\le k=|J| \le r-1}}
\rho_J(u(I))
\prod_{i\in I\cap J}(t_i-1)
(t_{r+1}-1)^{r-|I\cap J|-1}}\hfill\medskip\\
& & {\displaystyle 
\cdot \left\{\left( \prod_{i\in I}t_it_{r+1}^{-\sigma_J(u(I))}-1\right)
f_J(u(I))
+(t_{r+1}-1)h_J(u(I))\right\}}\hfill\medskip\\
& & {\displaystyle +\rho_{I_r}(u(I))\prod_{i\in I}(t_i-1)
(t_{r+1}-1)^{r-s}
h(u(I))}\hfill
\end{matrix}
\end{equation}
where $h_J$ and $h$ are the same as in Lemma \ref{lm:form}.
In (\ref{eq:tor3}),
we note that 
$$0\le |I\cap J|=|J|-|J\setminus I|\le |I|=s,$$
and
$|I\cap J|=s$ if and only if $J\supset I$.
By the fact, and (\ref{eq:tor1}), (\ref{eq:tor2}) and (\ref{eq:tor3}),
we have the result.
\qed

\bigskip

The Reidemeister torsions of $M_{r(I)}$ are the following:

\begin{lm}\label{lm:R-tor2}
Suppose that $r \ge 2$ and $1\le |I| \le r-1$.
\begin{enumerate}
\item[(1)]
If $I=\{x\}$\ $(s=1)$, then we have:
\begin{equation*}
\begin{matrix}
\tau(M_{u(I)}) & \doteq &
{\displaystyle 
\left\{
{\mit \Delta}_{K_x}(t_x)
+(t_x-1)^2\sum_{
{\scriptstyle J\ni x}
\atop
{\scriptstyle 2\le |J| \le r}}
\eta_J(u(I))F_J(I)\right\}(t_x-1)^{-1}}.\hfill
\end{matrix}
\end{equation*}

\item[(2)]
If $2\le |I| \le r-1$, then we have:
\begin{equation*}
\begin{matrix}
\tau(M_{u(I)}) & \doteq &
{\displaystyle \prod_{i\in I}(t_i-1)
\sum_{
{\scriptstyle J\supset I}
\atop
{\scriptstyle 2\le |J| \le r}}
\eta_J(u(I))F_J(I)}.\hfill
\end{matrix}
\end{equation*}

\end{enumerate}
\end{lm}

\noindent
{\bf Proof}\ 
By (\ref{eq:relation1}) and (\ref{eq:relation3}), we have 
$[m_{r+1}]=t_{r+1}=1$ and
$[l_{r+1}]=\prod_{i\in I}t_i$ in $H_1(M_{u(I)})$.
Hence by Lemma \ref{lm:surgery}, we have
$$\tau(M_{u(I)})\doteq
\tau(M_{u(\overline{I})})|_{t_{r+1}=1}
\left(\prod_{i\in I}t_i-1
\right)^{-1}.$$
By combining with Lemma \ref{lm:R-tor1},
we have the result.
\qed

\bigskip

For $I\subset I_r$ with $1\le |I|\le r$, and any 
$u(I)=(u_i)_{i\in I_r\setminus I}$,
we set
$$S_{u(I)}^{\mathrm{even}}
=\sum_{
{\scriptstyle J\supset I, |J\setminus I|\mathrm{: even}}
\atop
{\scriptstyle 2\le |J| \le r}}
\eta_J(u(I))F_J(I),$$
$$S_{u(I)}^{\mathrm{odd}}
=\sum_{
{\scriptstyle J\supset I, |J\setminus I|\mathrm{: odd}}
\atop
{\scriptstyle 2\le |J| \le r}}
\eta_J(u(I))F_J(I),$$
and
$-u(I)$ is obtained from $u(I)$ by replacing $u_i$ into $-u_i$
for all $i\in I_r\setminus I$.

\begin{lm}\label{lm:minus}
$$S_{-u(I)}^{\mathrm{even}}
=S_{u(I)}^{\mathrm{even}}\quad
\mbox{and}\quad
S_{-u(I)}^{\mathrm{odd}}
=-S_{u(I)}^{\mathrm{odd}}.$$
\end{lm}

\noindent
{\bf Proof}\ 
For $J\supset I$,
since $k_J(u(I))+k_J(-u(I))=|J\setminus I|$,
we have
$$\eta_J(-u(I))=(-1)^{|J\setminus I|}\eta_J(u(I)),$$
and the results.
\qed

\bigskip

\noindent
{\bf Proof of Theorem \ref{th:MT1}}\ 
Suppose that $L$ is component-preservingly amphicheiral.
Then $M_{u(I)}$ is homeomorphic to $M_{-u(I)}$.

\medskip

\noindent
(1)\ By Lemma \ref{lm:R-tor2} (1) and Lemma \ref{lm:minus}, we have
$$\tau(M_{u(I)}) \doteq 
\left\{
{\mit \Delta}_{K_x}(t_x)
+(t_x-1)^2
\left( S_{u(I)}^{\mathrm{even}}+S_{u(I)}^{\mathrm{odd}}\right)
\right\}(t_x-1)^{-1}$$
and
$$\tau(M_{-u(I)}) \doteq 
\left\{
{\mit \Delta}_{K_x}(t_x)
+(t_x-1)^2
\left( S_{u(I)}^{\mathrm{even}}-S_{u(I)}^{\mathrm{odd}}\right)
\right\}(t_x-1)^{-1}.$$
Since $\tau(M_{u(I)}) \doteq \tau(M_{-u(I)})$ and
${\mit \Delta}_{K_x}(1)=1\ne 0$,
we have $S_{u(I)}^{\mathrm{odd}}=0$.

\medskip

\noindent
(2)\ By Lemma \ref{lm:R-tor2} (2) and Lemma \ref{lm:minus}, we have 
$$\tau(M_{u(I)}) \doteq 
\prod_{i\in I}(t_i-1)
\left( S_{u(I)}^{\mathrm{even}}+S_{u(I)}^{\mathrm{odd}}\right)$$
and
$$\tau(M_{-u(I)}) \doteq 
\prod_{i\in I}(t_i-1)
\left( S_{u(I)}^{\mathrm{even}}-S_{u(I)}^{\mathrm{odd}}\right).$$
Since $\tau(M_{u(I)}) \doteq \tau(M_{-u(I)})$,
we have $S_{u(I)}^{\mathrm{even}}=0$
or $S_{u(I)}^{\mathrm{odd}}=0$.
\qed

\bigskip

\noindent
{\bf Proof of Corollary \ref{co:F=0}}\ 
(1)\ We set $r=2r'$
where $r'\in \mathbb{Z}$ and $r'\ge 1$.
We prove by induction on $r$.

\medskip

Suppose $r=2$ ($r'=1$).
By Theorem \ref{th:MT1} (1), 
we have $S_{u(I)}^{\mathrm{odd}}=\pm F(I)=0$.

\medskip

Suppose $r' \ge 2$.
By the assumption of induction,
we have $F_J(I)=0$ for every $J$ such that 
$|J|$ is even, and $2\le |J|\le r-2$.
By Theorem \ref{th:MT1} (1), we have 
$S_{u(I)}^{\mathrm{odd}}=\pm F(I)=0$.

\medskip

\noindent
(2)\ Since $S_{u(I)}^{\mathrm{even}}=\pm {\mit \Delta}_{L_I}\ne 0$,
we have $S_{u(I)}^{\mathrm{odd}}=\pm F(I)=0$ 
by Theorem \ref{th:MT1} (2).
The equation $F(I)=0$ holds if and only if
$f$ is divisible by $t_i-1$ from Lemma \ref{lm:dual}.

\medskip

\noindent
(3)\ Since $S_{u(I)}^{\mathrm{odd}}=\pm F(I)\ne 0$,
we have $S_{u(I)}^{\mathrm{even}}=\pm {\mit \Delta}_{L_I}=0$ 
by Theorem \ref{th:MT1} (2).
\qed

\begin{re}
{\rm
By Corollary \ref{co:F=0} (3) and Lemma \ref{lm:sub} (1),
for an algebraically split component-preservingly amphicheiral link $L$,
if we can add one component $K'$ to $L$ satisfying that
$L'=L\cup K'$ is also an algebraically split component-preservingly amphicheiral link 
such that ${\mit \Delta}_{L'}$ is not divisible by $(t'-1)^2$
where $t'$ corresponds to a meridian of $K'$,
then we have ${\mit \Delta}_L=0$.
We hope that it is possible for the case that
$L$ is an algebraically split component-preservingly amphicheiral link
with even components.
If it is true, then Conjecture \ref{cj:cj1} is affirmative.
However, if $L$ is the Borromean rings ($3$-component link),
then there does not exist a knot like $K'$.}
\end{re}

\noindent
{\bf Proof of Corollary \ref{co:2-comp1}}\ 
We take $I=\{2\}$.
By Corollary \ref{co:F=0} (2) (or Corollary \ref{co:F=0} (1)),
$f(t_1, t_2)$ is divisible by $t_1-1$.
We can argue similarly for the case $I=\{1\}$.
Therefore we have the result.
\qed

\section{Proof of Theorem \ref{th:MT2}
and Corollary \ref{co:2-comp2}}\label{sec:MT2}

We prove Theorem \ref{th:MT2} by a slightly generalized argument
of Hartley \cite[Theorem 2.1]{Ha}.

\medskip

\noindent
{\bf Proof of Theorem \ref{th:MT2}}\ 
Suppose that $L$ is oriented.
We span a Seifert surface $F$ corresponding the orientation.
We set a Seifert matrix from $F$ as $S$.
We can compute the one variable Alexander polynomial of $L$ from $S$ as:
$$(t-1){\mit \Delta}_L(t, \ldots, t)=\det(tS-S^T)$$
where $S^T$ is the transposed matrix of $S$.
Let $\varphi$ be an orientation-reversing homeomorphism of $S^3$.
Since $L$ is $(\varepsilon, \ldots, \varepsilon)$-amphicheiral, 
$\varphi(F)$ is still a Seifert surface of $L$, and 
the corresponding Seifert matrix changed into $-S$.
Since the S-equivalences do not change the one variable Alexander polynomial,
we have
$$(t-1){\mit \Delta}_L(t, \ldots, t)=\det(-tS+S^T).$$
Since the size of $S$ is odd, we have ${\mit \Delta}_L(t, \ldots, t)=0$.
We can argue similarly if $K_i$ is changed into $\eta_iK_i$.
Therefore we have the result.
\qed

\bigskip

\noindent
{\bf Proof of Corollary \ref{co:2-comp2}}\ 
By Corollary \ref{co:2-comp1} and Theorem \ref{th:MT2},
we have the result.
\qed

\section{Supporting examples}\label{sec:rem}
We raise examples which support Conjecture \ref{cj:cj1}.

\begin{ex}\label{ex:ex4}
{\rm
(1)\ In Figure 1, let $M_{\lambda}$ be 
the $\lambda$-component Milnor link \cite{Mi1} where $\lambda \ge 3$.
In particular, $M_3$ is the Borromean rings.
The Alexander polynomial of $M_{\lambda}$ is
${\mit \Delta}_{M_3}(t_1, t_2, t_3)=(t_1-1)(t_2-1)(t_3-1)$ and 
${\mit \Delta}_{M_{\lambda}}(t_1, \ldots, t_{\lambda})=0$\ 
$(\lambda \ge 4)$.
The Borromean rings $M_3$ is $(+, +, +)$-amphicheiral,
but it is not $(-, -, -)$-amphicheiral (cf.\ \cite{JLWW, KK}).

\begin{figure}[htbp]
\begin{center}
\includegraphics[scale=0.7]{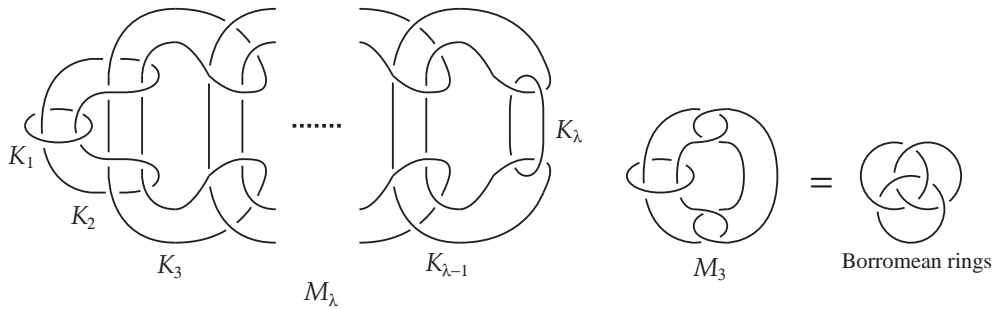} 
\label{fig:Milnor}
\caption{$\lambda$-component Milnor link $M_{\lambda}$}
\end{center}
\end{figure}

\noindent
(2)\ In Figure 2, let $C(2a_1, 2b_1, \ldots, 2a_n)$ be 
a $2$-component $2$-bridge link where the number in a rectangle
implies the number of half twists.
A $2$-component amphicheiral $2$-bridge link is not
algebraically split (see \cite{Kd1}), and the Alexander polynomial
of every non-trivial $2$-bridge link is not zero.
As a special case, 
the Alexander polynomial of $L=C(2a, 2b, -2a)$\ $(a\ne 0,\ b\ne 0)$ 
($C(2, \pm 2, -2)$ is the Whitehead link) is
$${\mit \Delta}_L(t_1, t_2)
=b(t_1-1)(t_2-1)\left\{ \frac{(t_1t_2)^a-1}{t_1t_2-1}\right\}^2$$
by Kanenobu's formula \cite{Kn}.
We can see that $L$ is not amphicheiral by Corollary \ref{co:2-comp1}.

\begin{figure}[htbp]
\begin{center}
\includegraphics[scale=0.6]{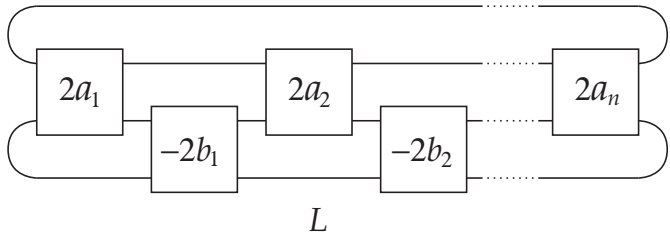} 
\label{fig:2-br}
\caption{$2$-bridge link $L=C(2a_1, 2b_1, \ldots, 2a_n)$}
\end{center}
\end{figure}

\noindent
(3)\ For links with up to $11$ crossings,
we use slightly modified notations in a web site 
maintaied by D.~Bar-Natan and S.~Morrison \cite{BM}.
In the class, only $10_{n36}^2$ and $10_{n107}^4$ are 
algebraically split component-preservingly amphicheiral links
with even components.
Moreover they are 
algebraically split component-preservingly $(+)$-amphicheiral links.
We can confirm that the Alexander polynomials of them are $0$
by direct computations or \cite[Theorem 1.3]{KK}.
We also remark that the condition ``component-preservingly" is needed.
$10_{n59}^2$ and $11_{n247}^2$
are algebraically split amphicheiral links with even components
which are not component-preservingly amphicheiral (cf.\ \cite{Kd3}).
The Alexander polynomials of them are
\begin{eqnarray*}
{\mit \Delta}_{10_{n59}^2}(t_1, t_2) & \doteq &
(t_1-1)(t_2-1)(t_1-t_2)(t_1t_2-1)
\medskip\\
{\mit \Delta}_{11_{n247}^2}(t_1, t_2) & = &
0.
\end{eqnarray*}
We note that
${\mit \Delta}_{10_{n59}^2}(t_1, t_2)$ 
satisfies the condition
$${\mit \Delta}_{10_{n59}^2}(t, t)
={\mit \Delta}_{10_{n59}^2}(t, t^{-1})=0$$
in Theorem \ref{th:MT2}, and both
$10_{n59}^2$ and $11_{n247}^2$ are $(\pm, \pm; (1\ 2))$-amphicheiral
where $(1\ 2)$ is the nontrivial permutation of $\{1, 2\}$.

\begin{figure}[htbp]
\begin{center}
\includegraphics[scale=0.6]{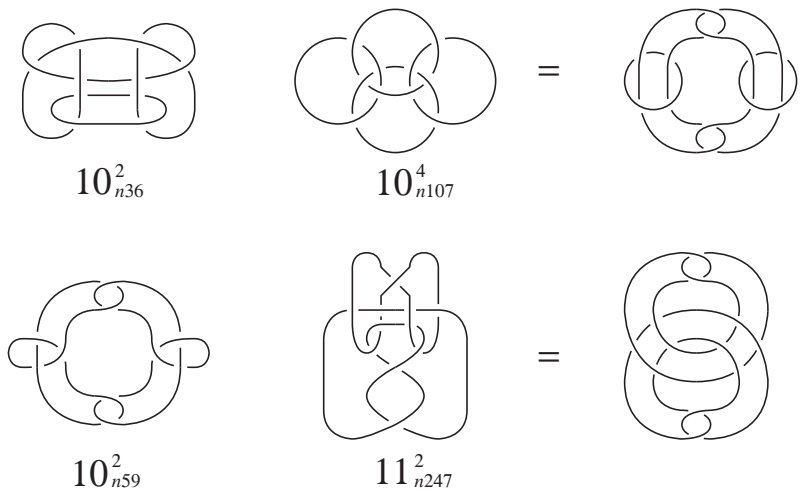} 
\label{fig:11cr}
\caption{Examples of prime links with up to $11$ crossings}
\end{center}
\end{figure}}
\end{ex}

{\noindent {\bf Acknowledgements}}\ 
The author would like to express gratitude to
Akio Kawauchi, Tsuyoshi Sakai and Shicheng Wang for giving him useful advices.

{\footnotesize
\par
Teruhisa KADOKAMI\par
Department of Mathematics,\par
East China Normal University,\par
Dongchuan-lu 500, Shanghai, 200241, China \par
{\tt mshj@math.ecnu.edu.cn}\par
{\tt kadokami2007@yahoo.co.jp}\par
}
\end{document}